\documentclass[intlimits]{amsart}
\usepackage{epsf,amscd}

\newcommand{\h}{Hoch\-schild\ }
\newcommand{\hc}{\h com\-plex}
\newcommand{\coh}{co\-ho\-mol\-o\-gy}

\newcommand{\ass}{\mathcal{A}ssoc}
\newcommand{\Aut}{\operatorname{Aut}}
\newcommand{\comm}{\mathcal{C}omm}
\newcommand{\End}[1]{\mathcal{E}nd\,_{#1}}
\newcommand{\G}{\mathcal{G}}
\newcommand{\Hom}{\operatorname{Hom}}
\newcommand{\id}{\operatorname{id}}
\newcommand{\im}{\operatorname{Im}}
\newcommand{\II}{\mathcal{I}}
\newcommand{\In}{\operatorname{In}}

\newcommand{\Ker}{\operatorname{Ker}}
\newcommand{\lie}{\mathcal{L}ie}
\newcommand{\Map}{\operatorname{Map}}
\newcommand{\met}{\operatorname{met}}
\newcommand{\MM}{\mathcal{M}}
\newcommand{\Mor}{\operatorname{Mor}}
\newcommand{\OO}{\mathcal{O}}

\newcommand{\PP}{\mathcal{P}}

\newcommand{\del}{\partial}
\newcommand{\nq}{{\mathbb{Q}}}
\newcommand{\nr}{{\mathbb{R}}}
\newcommand{\nz}{{\mathbb{Z}}}

\newtheorem{thm}{Theorem}[section]
\newtheorem{prop}[thm]{Proposition}
\newtheorem{crl}[thm]{Corollary}

\theoremstyle{definition}

\newtheorem{df}{Definition}[section]
\newtheorem{ex}{Example}[section]
\newtheorem{xca}{Exercise}

\theoremstyle{remark}

\newtheorem{rem}{Remark}%[section]

\newenvironment{ack}{\small \trivlist \item[\hskip \labelsep{\it
Acknowledgments}.]}{\endtrivlist}
\newenvironment{disc}{\small \trivlist \item[\hskip \labelsep{\it
Disclaimer}.]}{\endtrivlist}

\newcommand{\abs}[1]{\lvert#1\rvert}

\begin{document}

\title{Notes on Universal Algebra}

\author{Alexander A. Voronov} \address {School of Mathematics\\ 
  University of Minnesota\\Minneapolis, MN 55455}
\email{voronov@math.umn.edu}

\dedicatory{Dedicated to Dennis Sullivan on the occasion of his sixtieth
birthday.}
\date{October 25, 2001}
\thanks{The author was supported in part by NSF grant DMS-9971434 and the James~H. Simons Math-Physics Endowment}
\begin{abstract}
These are notes of a mini-course given at Dennisfest in June 2001. The
goal of these notes is to give a self-contained survey of deformation
quantization, operad theory, and graph homology. Some new results
related to ``String Topology'' and cacti are announced in
Section~\ref{cacti}.
\end{abstract}

\maketitle

\section*{Introduction}

Either due to the influence of string theory or just because this is
what the face of mathematics was supposed to look like before the
beginning of the third millennium, graphs have recently stepped
forward and overwhelmed many areas of mathematics, including universal
algebra. The use of graphs is similar to the Feynman diagram technique
in physics; the amazing thing is that its applications to pure
mathematics are extremely powerful. In these lectures we are going to
discuss three topics: deformation quantization, operad theory, and
graph homology, in which significant progress has been made with the
help of graphs as a unifying pattern:-) If you are in any way offended
with my choice of a title of this paper, please, see the disclaimer
below.

\begin{ack}
  I would like to thank Dennis Sullivan and the organizers of
  Dennisfest for inviting me to give this mini-course. These notes are
  based on numerous ideas of Maxim Kontsevich, which have largely
  influenced the subject. I am also grateful to Giovanni Felder, Kolya
  Ivanov, John Jones, Tom Leinster, Bob Penner, Jim Stasheff, Dennis
  Sullivan, and the lively Dennisfest and the University of Minnesota
  audiences for many helpful suggestions, incorporated in this version
  of the notes.  I would also like to thank Kyoto University for its
  hospitality in July 2001, when the paper was written up.
\end{ack}
%\end{sloppypar}

\begin{disc}
  The characters and figures portrayed and the titles and notions used
  herein are fictitious and any resemblance to the names, character,
  or history of any person, except Dennis Sullivan, is coincidental
  and unintentional.
\end{disc}

\section{Graphs and formal algebraic quantization}

One of the oldest open problems solved with the essential help of
graphs is perhaps the problem of deformation quantization.

\subsection{Deformations of associative algebras}

Let us start with a review of deformation theory of associative
algebras. Let $A$ be an associative algebra over a field $k$ of
characteristic zero.

\begin{df}
A \emph{formal deformation} of $A$ is a $k[[t]]$-bilinear
multiplication law $m_t: A[[t]] \otimes_{k[[t]]} A[[t]] \to A[[t]]$ on
the space $A[[t]]$ of formal power series in a variable $t$ with
coefficients in $A$, satisfying the following properties:
\[
m_t (a,b) = a \cdot b + m_1(a,b) t + m_2 (a,b) t^2 + \dotsb \qquad
\text{for $a, b \in A$,}
\]
where $a \cdot b$ is the original multiplication on $A$, and $m_t$ is
associative, which is equivalent to the equation
\[
m_t(m_t(a,b),c) = m_t (a, m_t(b,c)) \qquad \text{for $a,b,c \in A$.}
\]
\end{df}

\begin{rem}
\label{dq}
Note that for a formal deformation $m_t$ of a \emph{commutative}
algebra $A$, the bracket defined by the first-order part of the
commutator,
\begin{eqnarray*}
\{a,b\} & := & \frac{1}{2t} (m_t(a,b) - m_t(b,a)) \mod t,\\ & = &
  \frac{1}{2} (m_1(a,b) - m_1 (b,a))
\end{eqnarray*}
defines the structure of a Poisson algebra on $A$. (All the identities
of a Poisson algebra follow from the fact that $A[[t]]$ with the
product $m_t(a,b)$ and the bracket $\frac{1}{2} (m_t(a,b) - m_t(b,a))$
is a noncommutative Poisson algebra in an obvious sense, e.g., see
\cite{fgv}.)  In physical terms, one can regard $t$ as a quantum
parameter, such as the Planck constant, the Poisson algebra $A$ as the
\emph{quasi-classical limit} of the associative algebra $A[[t]]$, and
the algebra $A[[t]]$ as a
\emph{deformation quantization} of the Poisson algebra $A$.  The
\emph{deformation quantization problem} is the inverse problem: given
a Poisson algebra $A$, find a formal deformation returning the
original Poisson algebra structure on $A$ in the quasi-classical
limit.
\end{rem}

The main tool in studying deformation theory is the \emph{\hc}
\[
0 \to C^0(A,A) \xrightarrow{d} \dotsb \xrightarrow{d} C^n(A,A)
\xrightarrow{d} C^{n+1} (A,A) \xrightarrow{d} \dotsb,
\]
where $C^n(A,A) := \Hom(A^{\otimes n}, A)$ is the space of\emph{ \h
$n$-cochains}, \emph{i.e}., the $n$-linear maps $f(a_1, \dots, a_n)$
on $A$ with values in $A$, and the \emph{differential} $d$, $d^2 = 0$,
is defined as
\begin{equation*}
%\label{Hd}
\begin{split}
(d f) (a_1, \dots, a_{n+1}) 
:= \; &
\; a_1 f(a_2, \dots, a_{n+1}) \\
 &  + \sum_{i=1}^n (-1)^i f (a_1, \dots , a_{i-1}, a_i a_{i+1},
a_{i+2}, \dots , a_{n+1})\\
 & -  (-1)^{n} f (a_1, \dots, a_n) a_{n+1},\\
\end{split}
\end{equation*}
for $f \in C^n (A,A)$, $a_1, \dots, a_{n+1} \in A$. The \emph{\h \coh}
is then the \coh\
\[
H^\bullet(A,A) := \Ker d /\im d
\]
of this complex. The \hc\ admits a bracket
\[
[,]: C^m (A,A) \otimes C^n(A,A) \to C^{m+n-1} (A,A),
\]
called the \emph{Gerstenhaber bracket, or the G-bracket}. The formula
for this bracket, defined by M.~Gerstenhaber \cite{Gerst} is not
really inspirational to me, in spite of all those years I have spent
staring at it. A conceptual definition, due to J.~Stasheff, is based
on the following idea. The matter is that the \hc\ may be identified
with the space of graded derivations of the tensor coalgebra: $T^c
(A[1]) := \bigoplus_{n \ge 0} A[1]^{\otimes n}$, and the bracket is
just the commutator of derivations \cite{jim:bracket}. Here $A[1]$
denotes the graded vector space whose only nonzero graded component is
$A$, placed in degree $-1$. In general, for a graded vector space $V =
\bigoplus_n V^n$, $V[k]$ denotes \emph{grading shift}, or what is
known to topologists as \emph{$k$-fold suspension}: it is a graded
vector space $V[k]$ whose component of degree $n$ is $V[k]^n :=
V^{k+n}$. Note that a derivation determined by a map $A[1]^{\otimes n}
\to A[1]$ has degree $n-1$, therefore, the bracket defines a
differential graded (DG) Lie algebra structure on the \hc\ $C^\bullet
(A,A)[1]$ with a shifted grading $\deg f = n-1$ for $f \in C^n(A,A)$.
Note that the (appropriately shifted) \h \coh\ $H^\bullet (A,A)[1]$
inherits the structure of a graded Lie algebra.

The importance of the bracket comes from the following tautological
fact, which however may be regarded as the cornerstone of deformation
theory.
\begin{prop}
A formal multiplication
\[
m_t (a,b) = m_0(a, b) + m_1(a,b) t + m_2 (a,b) t^2 + \dotsb, \qquad a, b
\in A,
\]
is associative, iff $[m_t,m_t] = 0$.
\end{prop}

\begin{proof}
We need a little formula for the G-bracket of $m_t$ with itself:
\[
[m_t,m_t](a,b,c) = 2 (m_t(m_t(a,b),c) - m_t(a,m_t(b,c))),
\]
which would have been obvious, if we had used Gerstenhaber's formula
to define the G-bracket. It is a good exercise to deduce it directly
from the definition of the G-bracket as the commutator of
derivations. The right-hand side of the formula contains the
associativity equation, and we are done.
\end{proof}

\begin{rem}
\label{m_0}
Because the original multiplication $m_0(a,b) := a\cdot b$ is
associative, the G-bracket square of it vanishes: $[m_0,m_0] =
0$. Therefore, the commutator with $m_0$ defines an inner differential
on the \hc. This differential is in fact the \h differential: another
exercise is to verify that $d f = [f,m_0]$.
\end{rem}

Classical deformation theory was about the following ``perturbative''
results.
\begin{crl}
\begin{enumerate}
\item
A formal multiplication
\[
m_t (a,b) = a\cdot b + m_1(a,b) t + m_2 (a,b) t^2 + \dotsb, \qquad a,
b \in A,
\]
is associative modulo $t^2$, iff $d m_1 = 0$. In this case $m_1$
defines a \h \coh\ class $m_1 \in H^2(A,A)$.

\item
Suppose that a formal multiplication as above is associative modulo
$t^2$. Then the existence of $m_2$ such that $m_t$ is associative
modulo $t^3$ is equivalent to the vanishing $[m_1,m_1] = 0$ in \h
\coh\ $H^\bullet(A,A)$.

\end{enumerate}
\end{crl}

\begin{proof}
Expand the equation $\frac{1}{2}[m_t,m_t] = 0$ in powers of $t$ and
collect terms by $t^n$ for $n=0$, 1, and 2. We will get the following.
\begin{eqnarray*}
t^0: & & \frac{1}{2} [m_0,m_0] = 0,\\
t^1: & & d m_1 = 0,\\
t^2: & & dm_2 + \frac{1}{2} [m_1,m_1] = 0.
\end{eqnarray*}
The first equation is always satisfied, because the original
multiplication is associative, see Remark~\ref{m_0}. The other two
equations explain both statements of the corollary.
\end{proof}

\subsection{Deformation quantization}

Deformation quantization usually refers to a specific deformation
quantization problem in a geometric/physical setting. The following
theorem, which solves the deformation quantization problem posed by
Bayen, Flato, Fr{\o}nsdal, Lichnerowicz, and Sternheimer \cite{bffls},
is a remarkable breakthrough in pure mathematics achieved by applying
ideas motivated by Feynman diagrams.

\begin{thm}[M.~Kontsevich \cite{kon:dq}]
Every Poisson manifold $(M, \{, \})$ may be \emph{deformation
quantized}, \emph{i.e}., there exists a formal deformation
quantization, see Remark~\ref{dq},
\[
f \star g := m_t (f,g) = fg + m_1(f,g)t + m_2 (f,g) t^2 + \dotsb,
\qquad f, g \in C^\infty(M),
\]
of the Poisson algebra $A = C^\infty (M)$ of smooth functions, so that
all the $m_i$'s are \emph{local}, that is, bidifferential operators on
$M$. According to our definition of deformation quantization, the star
product must be associative and also recover the Poisson algebra of
functions in the quasi-classical limit, \emph{i.e.,} $(m_1(f,g) -
m_1(g,f))/2 = \{f,g\}$.
\end{thm}

\begin{proof}
We will only consider the case of $M = \nr^d$ with an arbitrary
Poisson structure, where the situation is already highly
nontrivial. Globalization, which is done using a Fedosov-type
connection, see \cite{kon:dq,cat-felder-toma}, lies outside the main
theme of these notes: no pattern in it has to do with graphs.

First, we will sketch Kontsevich's original proof, giving an explicit
formula for the \emph{star product} $f \star g$:
\begin{equation}
\label{star}
f \star g := \sum_{n=0}^\infty \frac{t^n}{n!} \sum_{\Gamma \in
\G_{n,2}} W_\Gamma B_\Gamma(f,g).
\end{equation}
which is explained in the this paragraph. The interior summation runs
over the set $\G_{n,2}$ of directed graphs $\Gamma$ of a certain type
with vertices labeled $1, 2, \dots, n, \bar 1, \bar 2$. The set
$\G_{n,2}$ of graphs consists of the graphs satisfying the following
conditions. Each vertex of first type, \emph{i.e}., labeled 1, 2,
\dots, or $n$, has exactly two outgoing edges, labeled the first one
and the second, and there are no other edges whatsoever. No edge may
form a loop, \emph{i.e}., start and end at one and the same
vertex. $B_\Gamma(f,g)$ is a bidifferential operator defined by an
explicit formula \cite{kon:dq}, which we will describe using an
example.

\centerline{\epsfxsize=2in \epsfbox{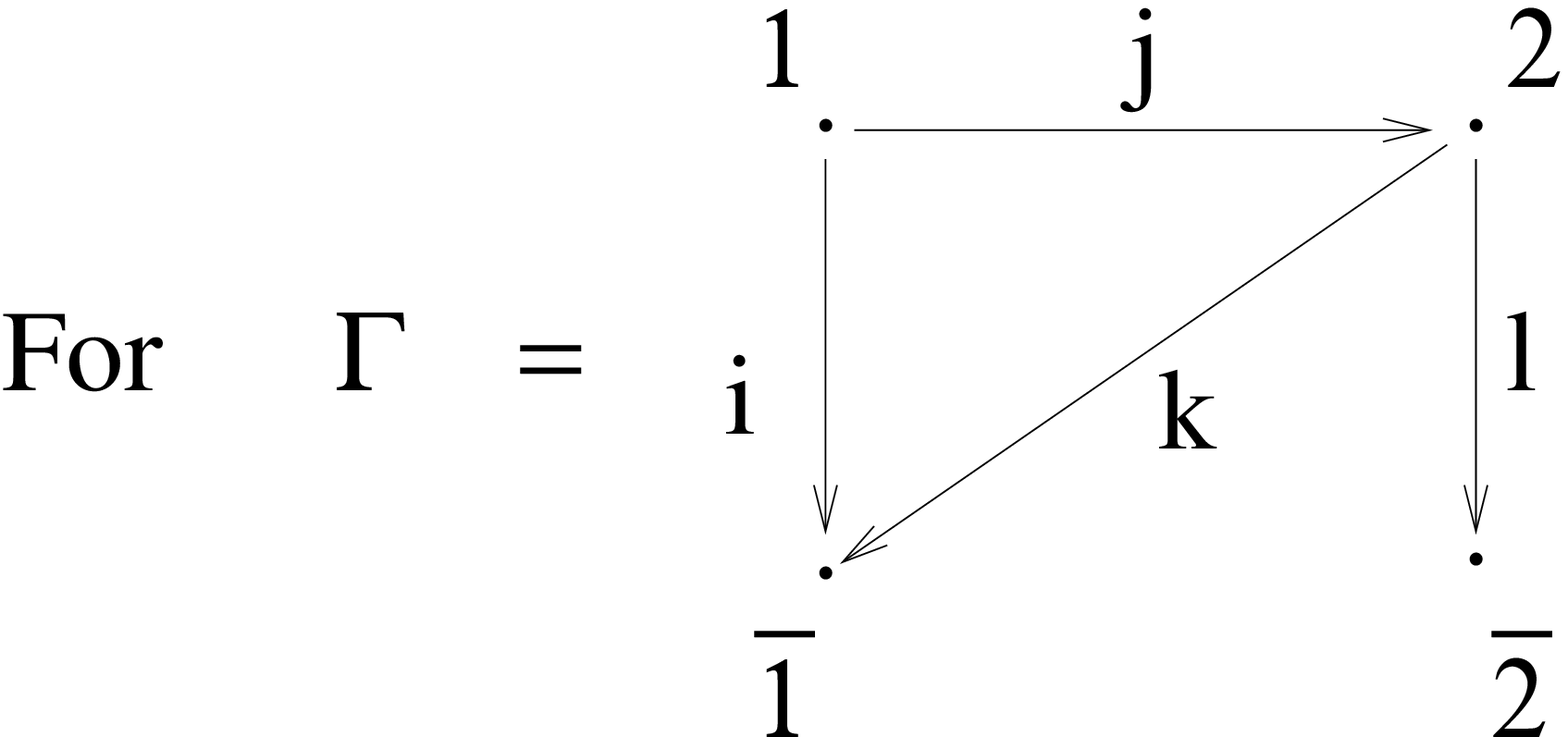}}
\smallskip
\noindent
the corresponding bidifferential operator will be
\[
B_\Gamma(f,g) := \alpha^{ij} \del_j(\alpha^{kl}) \del_k \del_i (f)
\del_l (g),
\]
where $\alpha^{ij}$ denotes the corresponding component of the Poisson
tensor in a fixed coordinate system $(x_1, \dots, x_d)$ in $\nr^d$ and
we assume summation over the repeating indices.  Finally, the
coefficient $W_\Gamma$ is given by another formula:
\[
W_\Gamma := \frac{1}{n! (2 \pi)^{2n}} \int_{C^+_{n,2}}
\bigwedge_{\text{edges $r \to s$ of $\Gamma$}} d \phi (z_r,z_s),
\]
where $C^+_{n,2}$ is the configuration space of $n$ distinct points
$z_1, \dots, z_n$ in the upper half-plane and two fixed points
$z_{\bar 1}=0 $ and $z_{\bar 2}=1$ on the real line; $\phi(z_r,z_s)$,
$r$ and $s$ running over $\{1,2, \dots, n, \bar 1, \bar 2\}$, is the
directed angle at $z_r$ between the hyperbolic line through $z_r$ and
$\infty$ and the hyperbolic line through $z_r$ and $z_s$. The order in
the wedge product is given by the lexicographic order of the vertices
$\{1, \dots, n\}$ and the orders of the set of edges going out of the
vertices. Kontsevich proves that this improper integral is absolutely
convergent.

The associativity of the star product may be verified explicitly, see
\cite{kon:dq}. However, I will sacrifice the rigor for the moral and
give a more conceptual, physical explanation of the associativity,
following A.~Cattaneo and G.~Felder \cite{cat-felder}. Cattaneo and
Felder define the star product as
\begin{equation}
\label{feynman}
(f \star g)(x) := \int_{\PP_x} f(X(0)) g(X(1)) e^{iS(X,\eta)/t} DX
D\eta.
\end{equation}
This is a Feynman integral over the infinite dimensional ``path''
space, which is the following space of fields $X$ and $\eta$ on the
upper half-plane $H$:
\begin{multline*}
\PP_x  := \{X: \bar H \to \nr^d, \eta \in \Omega^1 (\bar H) \otimes
\nr^d \; | \; X(\infty) = x \\
\text{ and $\eta$ vanishes on tangent vectors to the boundary}\}.
\end{multline*}
The function $S(X,\eta)$ is a certain action functional defining a
Poisson sigma model on $H$, see \cite{cat-felder}:
\[
S = \int_H (\eta_{\mu i} \del_\nu X^i + \frac{1}{2} \alpha^{ij}(X)
\eta_{\mu i} \eta_{\nu j} ) du^\mu du^\nu.
\]
A rigorous definition of the Feynman integral would be the very
formula \eqref{star}. However, physics takes the opposite viewpoint
and treats \eqref{star} as the saddle-point expansion of the integral
\eqref{feynman} in parameter $t$ obtained by formally applying the
rules by analogy with the finite-dimensional case. The advantage of
this approach is that the mystery of Kontsevich's formula \eqref{star}
is now replaced by the mystery of Equation \eqref{feynman}, which is
not so mysterious to a physicist, for whom it represents a standard
integral quantization formula. Another advantage is that it offers the
following explanation of the associativity.

Consider the integral
\[
\langle f, g , h  \rangle_p (x) := \int_{\PP_{x}} f(X(0)) g(X(1)) h(X(p))
e^{i S(X,\eta)/t} DX D\eta,
\]
where $p \in (1, \infty) \subset \nr \subset \bar H$ is a fixed point
on the real line between 1 and $\infty$ and $\PP_x$ is as above in
\eqref{feynman}. This integral is independent of the choice of this
point $p$, because the action $S$ is diffeomorphism invariant and,
roughly speaking, by integrating over all fields $X$ and $\eta$, we
take an average over all possible positions of $p$. Thus, the limits
of $\langle f,g,h \rangle_p$ as $p \to 1$ and $p \to \infty$ will be
the same. On the other hand, in the moduli space of configurations of
four points $0, 1, p$, and $\infty$ on the boundary of $H$, these
configurations will degenerate as follows:
\smallskip

\centerline{\epsfxsize=3in \epsfbox{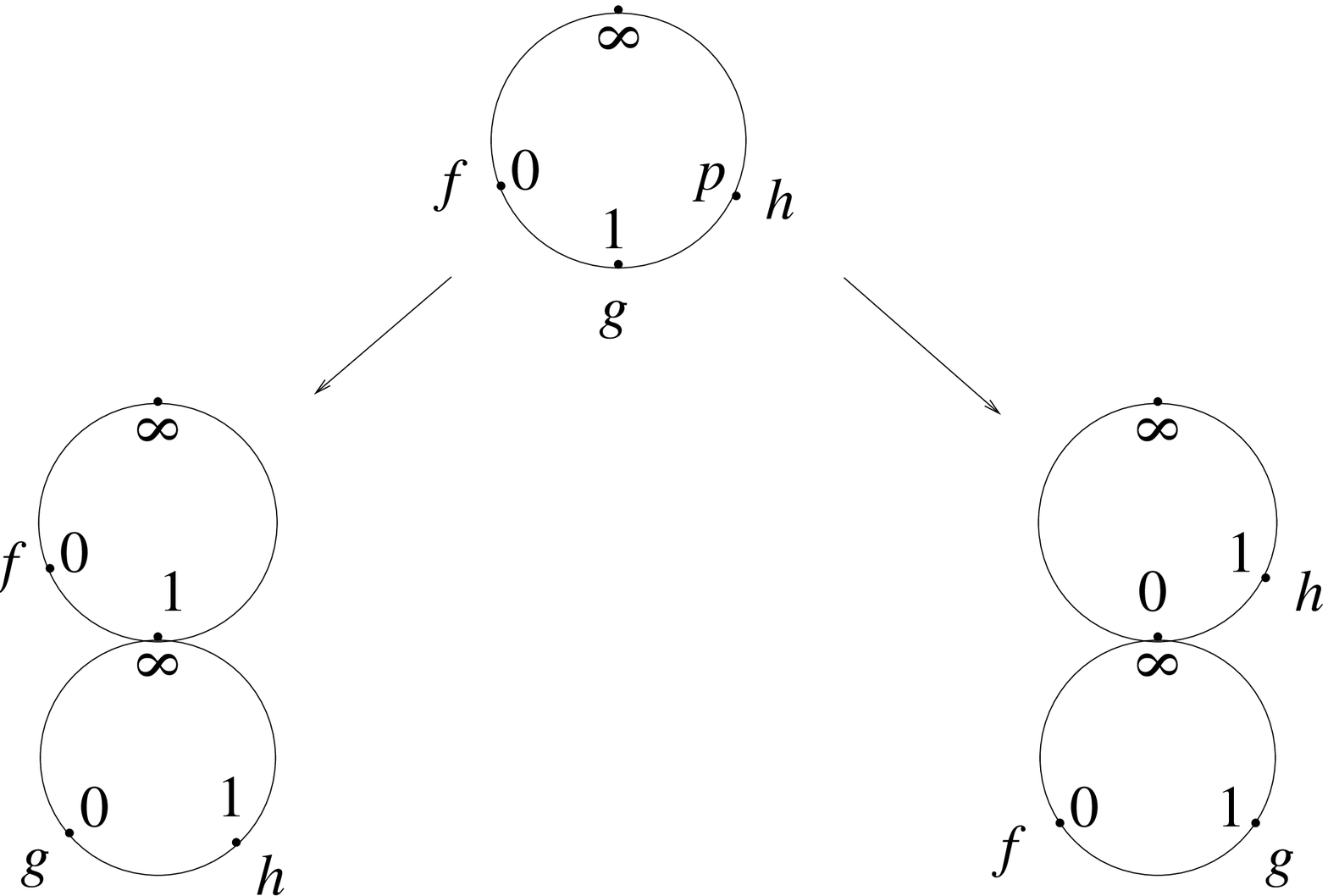}}
\smallskip
\noindent
which means that
\begin{eqnarray*}
\lim_{p \to 1} \langle f, g , h  \rangle_p  & = & f \star (g \star h),\\
\lim_{p \to \infty} \langle f, g , h  \rangle_p  & = & (f \star g) \star h,
\end{eqnarray*}
yielding the associativity.

In reality, things are more complicated than I have made them appear:
the Feynman diagram expansion involves gauge fixing and
renormalization, which is achieved by introducing ghosts and
antighosts (and what not) and using the BV formalism, see
\cite{cat-felder}. The behavior of the Feynman integral with respect
to the compactification of the configuration spaces is another issue
suppressed in the above. I hope to learn these things before the next
Dennisfest:-)

\end{proof}

\section{Trees and operads}

The combinatorics of trees essentially describes the combinatorics of
operads, while more general graphs would be related to PROP's and
modular operads. The main idea is that the trees form a free operad
(More precisely, the linear span of the set of labeled trees with a
certain differential forms the $A_\infty$ operad), and vice versa, any
free operad may be described via decorated trees: the tree remembers
how to get to its root by applying the operations put at its vertices,
given no identities, such as the associativity, whatsoever. Here we
will review these notions and related concepts of homotopy algebra.

Operads in general are spaces of operations with certain rules on how
to compose operations to get new ones. In this sense operads are
directly related to Lawvere's algebraic theories and represent true
objects of universal algebra. However, operads as such appeared in
topology in the works of J.~P. May, J.~M. Boardman and R.~M. Vogt as a
recognition tool for based multiple loop spaces. Stasheff earlier
described the first example of an operad, the associahedra, which
recognized based loop spaces, while about the same time, Gerstenhaber,
studying the algebra of the Hochschild complex, introduced the notion
of a composition algebra, which is nothing but the notion of an operad
of graded vector spaces.

\subsection{PROP's}

Unfortunately, the formal definition of an operad resembles in a way
the definition of a Deligne-Mumford stack --- it is very easy to drown
in the ocean of sheaves of groupoids over the site of schemes on your
first trip to the beach. It is much easier to define a PROP (=PROducts
and Permutations) and think of an operad as certain part of a PROP.

\begin{df}
  A \emph{PRO} is a symmetric monoidal (\emph{i.e}., tensor) category
  whose set of objects is identified with the set $\nz_+$ of
  nonnegative integers, where the tensor law on $\nz_+$ is given by
  addition. A \emph{PROP} is a PRO along with a right action of the
  permutation group $S_m$ and a left action of $S_n$ on Mor$(m,n)$ for
  each pair $m, n \ge 0$. These actions should be in a natural way
  compatible with the tensor structure and the composition of
  morphisms in the category that constitutes the PROP, see, for
  example, the founding fathers' sources, such as, J.~F. Adams' book
  \cite{Ad} or S.~Mac Lane's paper \cite{maclane}, or a postmodern view\\
  \texttt{http://www.math.umn.edu/$\sim$voronov/8390/lec2.pdf}\\
  for more detail.
\end{df}

Usually, PROP's are enriched over another symmetric monoidal category,
that is, the morphisms in the PROP are taken as objects of the other
symmetric monoidal category. This gives the notions of a PROP of sets,
vector spaces, complexes, topological spaces, manifolds, etc.
Examples of PROP's include the following. We will only specify the
morphisms, because the objects are already given by the definition.

\begin{ex}
  The \emph{endomorphism PROP} of a vector space $V$ has the space of
  morphisms $\Mor (m,n) = \Hom(V^{\otimes m}, V^{\otimes n})$. This is
  a PROP of vector spaces. The composition and tensor product of
  morphisms is given by same of linear maps.
\end{ex}

\begin{ex}
  The \emph{Segal PROP} is a PROP of infinite dimensional complex
  manifolds. The space of morphisms is defined as the moduli space
  $\PP_{m,n}$ of complex Riemann surfaces bounding $m+n$ labeled
  nonoverlapping holomorphic holes. The surfaces should be understood
  as compact smooth complex curves, not necessarily connected, along
  with $m+n$ biholomorphic maps of the closed unit disk to the
  surface. The more exact nonoverlapping condition is that the closed
  disks in the inputs do not intersect pairwise and the closed disks
  in the outputs do not intersect pairwise, however, an input and an
  output disk may have common boundary, but are still not allowed to
  intersect at an interior point. This technicality brings in the
  identity morphisms to the PROP, but does not create singular Riemann
  surfaces by composition. The moduli space means that we consider
  isomorphism classes of such objects. The composition of morphisms in
  this PROP is given by sewing the Riemann surfaces along the
  boundaries, using the equation $zw=1$ in the holomorphic parameters
  coming from the standard one on the unit disk. The tensor product of
  morphisms is the disjoint union. This PROP plays a crucial role in
  Conformal Field Theory, as we will see now.
\end{ex}

\subsection{Algebras over a PROP}

We need to define another important notion before we proceed.
\begin{df}
  We say that a vector space $V$ is an \emph{algebra over a PROP} $P$,
  if a morphism of PROP's from $P$ to the endomorphism PROP of $V$ is
  given. A morphism of PROP's should be a functor respecting the
  symmetric monoidal structure and the symmetric group actions, and
  also equal to the identity on the objects.
\end{df}

An algebra over a PROP could have been called a \emph{representation},
but since algebras over operads, which are similar objects, are
nothing but familiar types of algebras, it is more common to use the
term ``algebra''.

\begin{ex}
An example of an algebra over a PROP is a \emph{Conformal Field Theory
$($CFT$)$}, which may be defined as an algebra over the Segal
PROP. The fact that the functor respects compositions of morphisms
translates into the sewing axiom of CFT in the sense of
G.~Segal. Usually, one also asks for the functor to depend smoothly on
the point in the moduli space $\PP_{m,n}$. This definition of a CFT
describes only theories with a vanishing central charge. One needs to
extend the Segal PROP by a line bundle to cover the case of an
arbitrary charge, see \cite{huang:book}.
\end{ex}

\begin{ex}[Sullivan]
  Another example of an algebra over a PROP is a Lie bialgebra.
  Sullivan has shared with me a nice graph description of the
  corresponding PROP, see\\
  \texttt{http://www.math.umn.edu/$\sim$voronov/8390/lec4.pdf}.
\end{ex}

\subsection{Operads}

Now we are ready to deal with operads. Informally, an operad is the
part $\Mor(n,1)$, $n \ge 0$, of a PROP. Of course, given only the
collection of morphisms $\Mor(n,1)$, it is not clear how to compose
them. The idea is to take the union of a $m$ elements from $\Mor(n,1)$
to be able to compose them with an element of $\Mor(m,1)$. This leads
to cumbersome notation and ugly axioms, compared to those of a
PROP. However operads are in a sense more basic than the corresponding
PROP's: the difference is similar to the difference between Lie
algebras and the universal enveloping algebras.

\begin{df}
An \emph{operad} $\OO$ is a collection of sets (vector spaces,
complexes, topological spaces, manifolds, \dots, objects of a
symmetric monoidal category) $\OO (n)$, $n \ge 0$, with
\begin{enumerate}
\item A composition law:
\end{enumerate}
\noindent
\begin{equation*}
\gamma: \OO(m) \otimes \OO(n_1) \otimes \dotsb \otimes \OO(n_m) \to
\OO(n_1 + \dotsb + n_m).
\end{equation*}
\begin{enumerate}
\setcounter{enumi}{1}
\item A right action of the symmetric group $S_n$ on $\OO (n)$.
\item A unit $e \in \OO (1)$.
\end{enumerate}
such that the following properties are satisfied:
\begin{enumerate}
\item The composition is associative, \emph{i.e}., the following
diagram is commutative:
\end{enumerate}
\noindent
\[
\begin{CD}
\left\{
\begin{aligned}
\OO(l) & \otimes \OO(m_1) \otimes \dotsb \otimes \OO(m_l) \\
& \otimes \OO(n_{11}) \otimes \dotsb \otimes \OO(n_{l,n_l})
\end{aligned}
\right\}
@>{\id \otimes \gamma^l}>>
\OO(l) \otimes \OO(n_1) \otimes \dotsb \otimes \OO(n_l)\\
@V{\gamma \otimes \id}VV        @VV{\gamma}V , \\
\OO(m) \otimes \OO(n_{11}) \otimes \dotsb \otimes \OO(n_{m,n_m})
@>\gamma>> \OO(n)
\end{CD}
\]
\noindent
\begin{quote}
where $m = \sum_i m_i$, $n_i = \sum_j n_{ij}$, and $n = \sum_i n_i$.
\end{quote}
\noindent
\begin{enumerate}
\setcounter{enumi}{1}
\item The composition is equivariant with respect to the symmetric
  group actions: the groups $S_m$, $S_{n_1}$, \dots, $S_{n_m}$ act on
  the left-hand side and map naturally to $S_{n_1 + \dotsb + n_m}$,
  acting on the right-hand side.
\item The unit $e$ satisfies natural properties with respect to the
composition: $\gamma(e; f) \linebreak[0] = f$ and $\gamma(f; e, \dots,
e) = f$ for each $f \in \OO (k)$.
\end{enumerate}

The notion of a \emph{morphism of operads} is introduced naturally.
\end{df}

\begin{rem}
One can consider \emph{non-$\Sigma$ operads}, not assuming the action
of the symmetric groups. Not requiring the existence of a unit $e$, we
arrive at \emph{nonunital operads}. Do not mix this up with operads
with no $\OO(0)$, algebras over which (see next section) have no
unit. There are also good examples of operads having only $n \ge 2$
components $\OO(n)$.

An equivalent definition of an operad may be given in terms of
operations $f \circ_i g \linebreak[1] = \linebreak[0] \gamma (f; \id,
\dots, \linebreak[1] \id, g, \id , \dots, \id)$, $i = 1, \dots, m$,
for $f \in \OO(m), g \in \OO(n)$. Then the associativity condition
translates as $f \circ_i (g \circ_j h) = (f \circ_i g) \circ_{i+j-1}
h$ plus a natural symmetry condition for $(f \circ_i g) \circ_j h$,
when $g$ and $h$ ``fall into separate slots'' in $f$, see \emph{e.g}.,
\cite{ksv2}.
\end{rem}

%\begin{sloppypar}
\begin{ex}[The Riemann surface and the endomorphism operads]
$\PP(n)$ is the space of Riemann spheres with $n+1$ boundary
components, \emph{i.e}., $n$ inputs and 1 output.  Another example is
the \emph{endomorphism operad of a vector space $V$}: $\End{V}(n) =
\Hom(V^{\otimes n}, V)$, the space of $n$-linear mappings from $V$ to
$V$.
\end{ex}
%\end{sloppypar}

\subsection{Algebras over an operad}

\begin{df}
An \emph{algebra over an operad $\OO$} (in other terminology, a
\emph{representation of an operad}) is a morphism of operads $\OO \to
\End{V}$, that is, a collection of maps
\[
\OO(n) \to \End{V}(n) \qquad \text{for $ n \ge 0$}
\]
compatible with the symmetric group action, the unit elements, and the
compositions. If the operad $\OO$ is an operad of vector spaces, then
we would usually require the morphism $\OO \to \End{V}$ to be a
morphism of operads of vector spaces. Otherwise, we would think of
this morphism as a morphism of operads of sets. Sometimes, we may also
need a morphism to be continuous or respect differentials, or have
other compatibility conditions.
\end{df}

\subsubsection{The commutative operad}
The \emph{commutative operad} is the operad of $k$-vector spaces with
the $n$th component $\comm (n) = k$ for all $n \ge 0$. We assume that
the symmetric group acts trivially on $k$ and the compositions are
just the multiplication of elements in the ground field $k$. An
algebra over the commutative operad is nothing but a commutative
associative algebra with a unit, as we see from the following
exercise.

Another version of the commutative operad is $\comm(n) = \text{point}$
for all $n \ge 0$. This is an operad of sets. An algebra over it is
also the same as a commutative associative unital algebra.

\begin{xca}
Show that the operad $\mathcal{T}op (n) \linebreak[1] = \linebreak[0]$
\{the set of diffeomorphism classes of Riemann spheres with $n$ input
holes and 1 output hole\} is isomorphic to the commutative operad of
sets.
\end{xca}

\begin{xca}
Prove that the structure of an algebra over the commutative operad
$\comm$ on a vector space is equivalent to the structure of a
commutative associative algebra with a unit.
\end{xca}

\subsubsection{The associative operad}
\label{ass:operad}

The \emph{associative operad} $\ass$ can be considered as a
one-dimensional analogue of the commutative operad
$\mathcal{T}op$. $\ass(n)$ is the set of equivalence classes of
connected planar binary (each vertex being of valence 3) trees that
have a root edge and $n$ leaves labeled by integers 1 through $n$:
\medskip

\centerline{\epsfysize=3cm \epsfbox{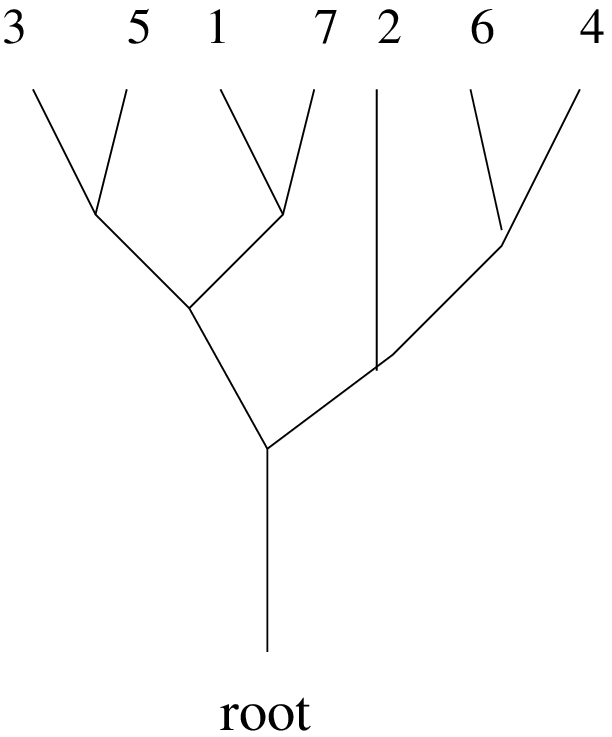}}

\noindent
If $n=1$, there is only one tree --- it has no vertices and only one
edge connecting a leaf and a root. If $n=0$, the only tree is the one
with no vertices and no leaves --- it only has a root. Unfortunately,
I have a problem sketching it: it probably exists only in the quantum
world.

Two trees are equivalent if they are related by a sequence of moves of
the kind
\medskip

\centerline{\epsfysize=2cm \epsfbox{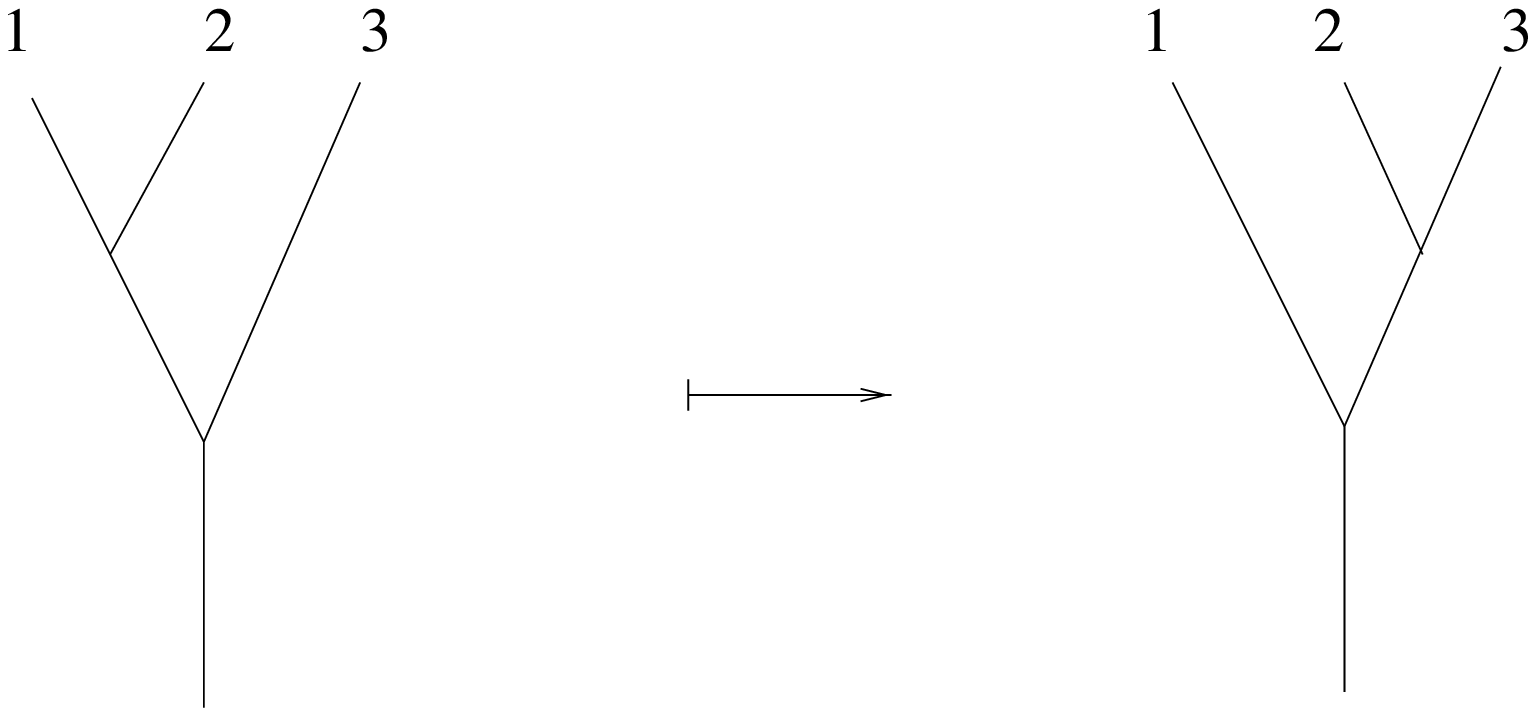}}
\smallskip

\noindent
performed over pairs of two adjacent vertices of a tree. The symmetric
group acts by relabeling the leaves, as usual. The composition is
obtained by grafting the roots of $m$ trees to the leaves of an
$m$-tree, no new vertices being created at the grafting points. Note
that this is similar to sewing Riemann surfaces and erasing the seam,
just as we did to define operad composition in that case. By
definition, grafting a 0-tree to a leaf just removes the leaf and, if
this operation creates a vertex of valence 2, we should erase the
vertex.

\begin{xca}
Prove that the structure of an algebra over the associative operad
$\ass$ on a vector space is equivalent to the structure of an
associative algebra with a unit.
\end{xca}

\subsubsection{The Lie operad}
\label{lie:operad}

The \emph{Lie operad} $\lie$ is another variation on the theme of a
tree operad. Consider the same planar binary trees as for the
associative operad, except that we do not include a 0-tree,
\emph{i.e}., the operad has only positive components $\lie(n)$, $n \ge
1$, and there are now two kinds of equivalence relations:
\medskip

\centerline{\epsfysize=2.5cm \epsfbox{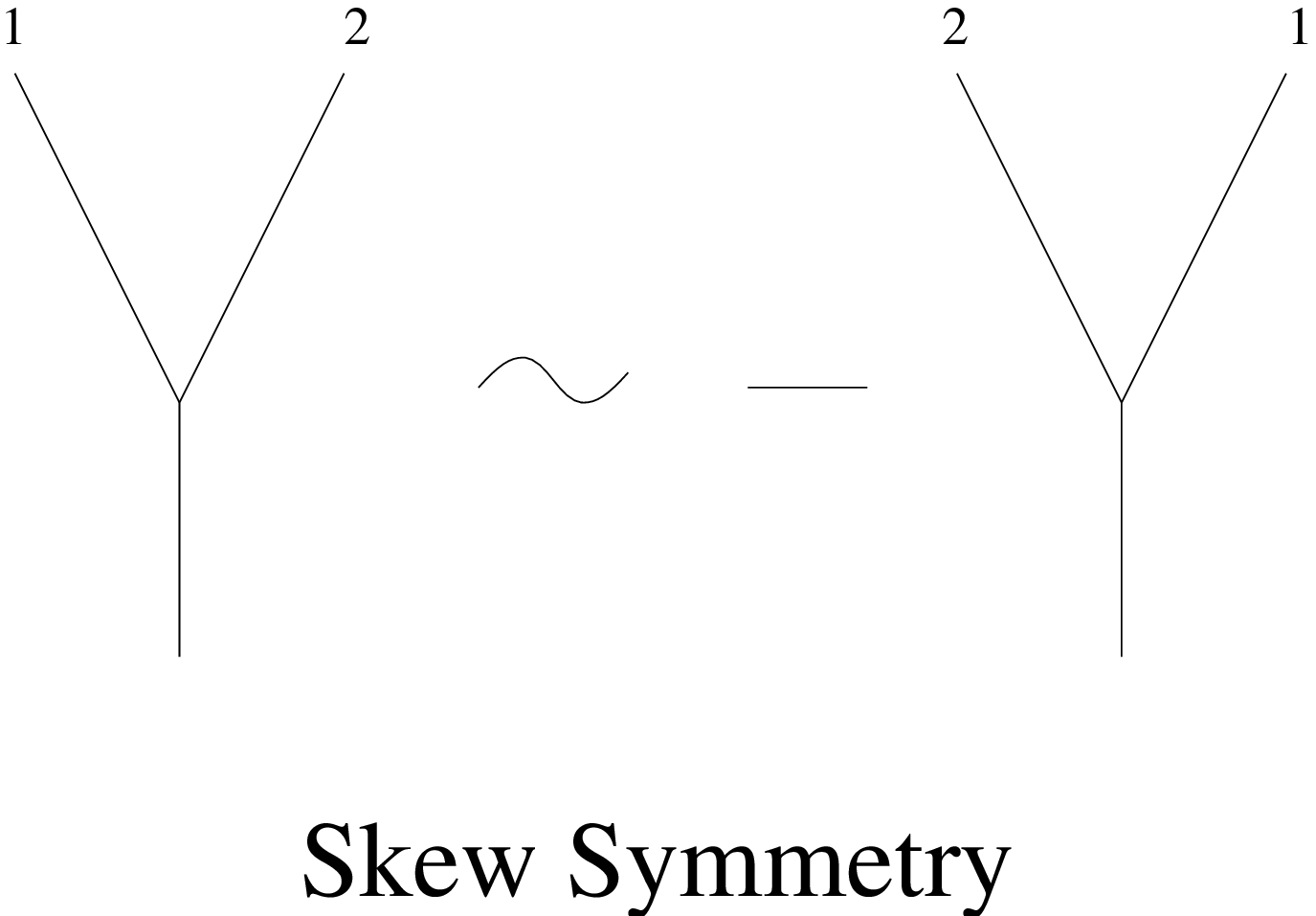}}
\smallskip

\noindent
and
\medskip

\centerline{\epsfysize=2.5cm \epsfbox{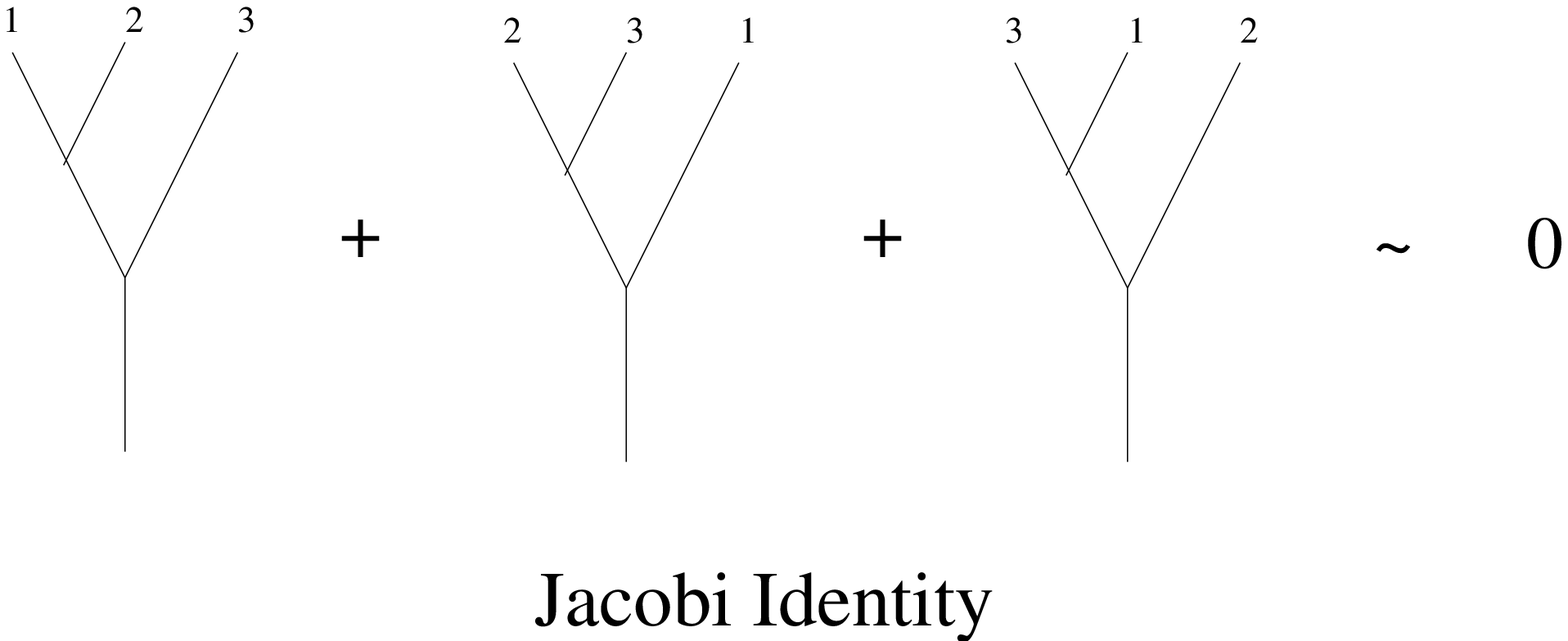}}
\smallskip

\noindent
Now that we have arithmetic operations in the equivalence relations,
we consider the Lie operad as an operad of vector spaces. We also
assume that the ground field is of a characteristic other than 2,
because otherwise we will arrive at the wrong definition of a Lie
algebra.

\begin{xca}
  Prove that the structure of an algebra over the Lie operad $\lie$ on
  a vector space over a field of a characteristic other than 2 is
  equivalent to the structure of a Lie algebra.
\end{xca}

\begin{xca}
Describe algebraically an algebra over the operad $\lie$, if we modify
it by including a 0-tree, whose composition with any other tree is
defined as (a) zero, (b) the one for the associative operad.
\end{xca}

\subsubsection{The Poisson operad}
\label{poisson:operad}

Recall that a \emph{Poisson algebra} is a vector space $V$ (over a
field of characteristic zero) with a unit element $e$, a dot product
$ab$, and a bracket $[a,b]$ defined, so that the dot product defines
the structure of commutative associative unital algebra, the bracket
defines the structure of a Lie algebra, and the bracket is a
derivation of the dot product:
\[
[a,bc] = [a,b] c + b [a,c] \qquad \text{for all $a$, $b$, and $c \in
V$}.
\]

\begin{xca}
Define the \emph{Poisson operad}, using a tree model similar to the
previous examples. Show that an algebra over it is nothing but a
Poisson algebra. [\emph{Hint}: Use two kinds of vertices, one for the
dot product and the other one for the bracket.]
\end{xca}

\subsubsection{The Riemann surface operad and vertex operator algebras}

Just for a change, let us return to the operad $\PP$ of Riemann
surfaces, more exactly, isomorphism classes of Riemann spheres with
holomorphic holes. What is an algebra over it? Since there are
infinitely many nonisomorphic pairs of pants, there are infinitely
many (at least) binary operations. In fact, we have an infinite
dimensional family of binary operations parameterized by classes of
pairs of pants. However modulo the unary operations, those which
correspond to cylinders, we have only one fundamental binary operation
corresponding to a fixed pair of pants. An algebra over this operad
$\PP$ is part of a CFT data at the tree level, the central charge
$c=0$. If we consider a holomorphic algebra over this operad, that is,
require that the defining mappings $\PP(n) \to \End{V}(n)$, where $V$
is a complex vector space, be holomorphic, then we get part of a
chiral CFT, or an object which may be called a \emph{vertex operator
  algebra $($VOA$)$}. This kind of object is not equivalent to what
people used to call a VOA, but according to Y.-Z. Huang's Theorem, a
true VOA is a holomorphic algebra over a ``partial pseudo-operad of
Riemann spheres with rescaling'', which is a version of $\PP$, where
the disks are allowed to overlap. The fundamental binary operation
$Y(a,z)b$ for $a, b \in V$ of a VOA is commonly chosen to be the one
corresponding to a pair of pants which is the Riemann sphere with a
standard holomorphic coordinate and three unit disks around the points
0, $z$, and $\infty$ (No doubt, these disks overlap badly, but we
shrink them on the figure to look better): \medskip

\centerline{\epsfysize=1.8cm \epsfbox{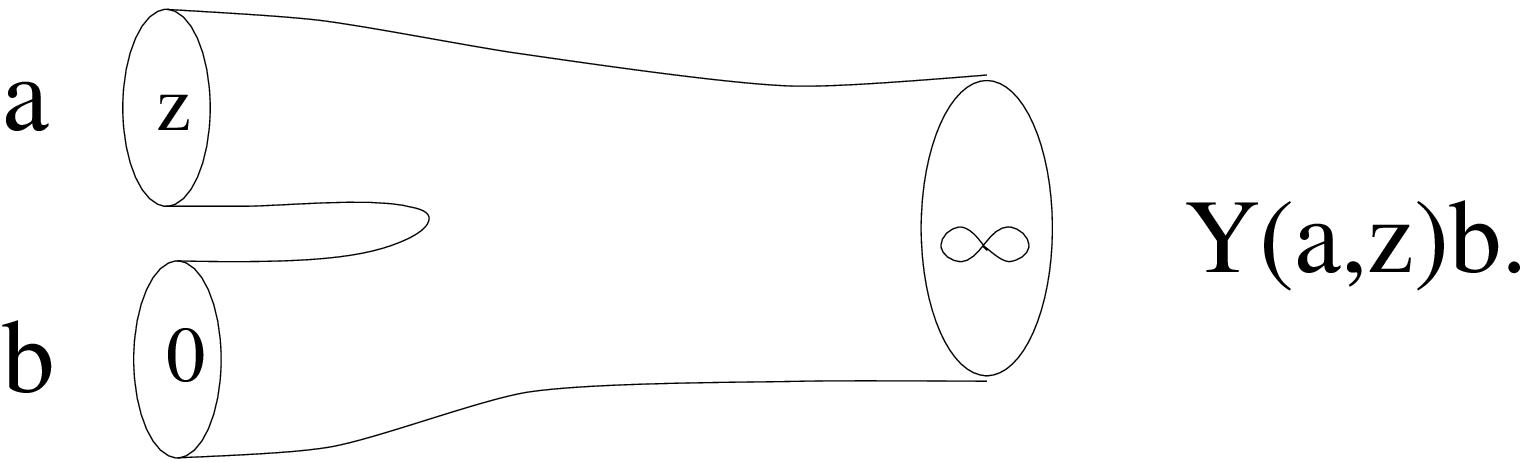}}
\smallskip

The famous associativity identity 
\[
Y(a,z-w)Y(b,-w) c = Y(Y(a,z)b,-w) c
\]
for vertex operator algebras comes from the following natural
isomorphism of the Riemann surfaces:
\medskip

\centerline{\epsfysize=2.5cm \epsfbox{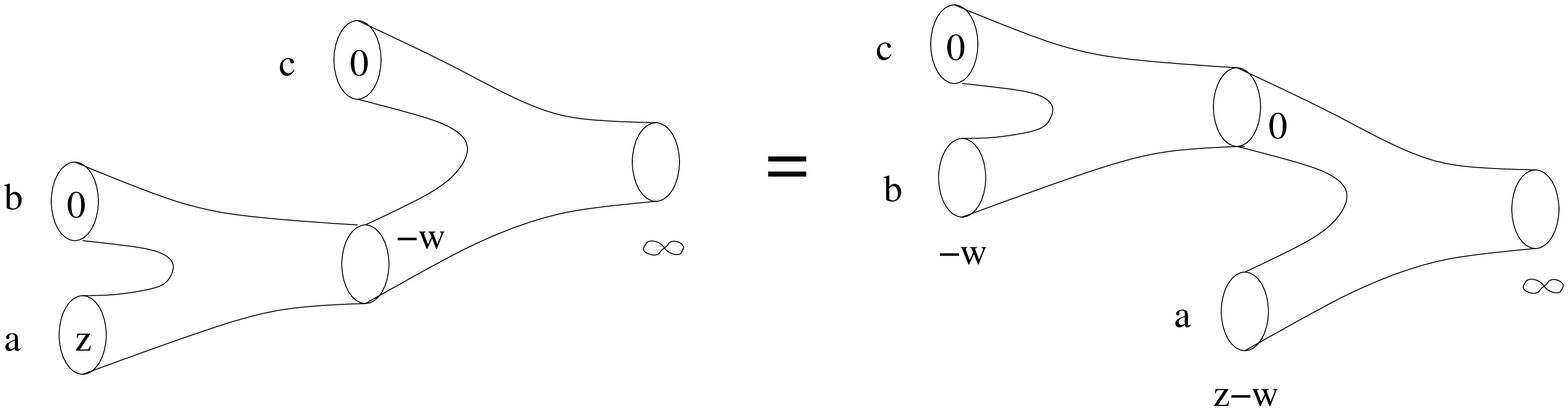}}
\smallskip

\subsection{Operads via generators and relations}

The tree operads that we looked at above, such as the associative and
the Lie operads, are actually operads defined by generators and
relations. Here is a way to define such operads in general. To fix
notation, assume throughout this section that we work with operads
$\OO(n)$, $ n \ge 1$, of vector spaces.

\begin{df}
An \emph{ideal} in an operad $\OO$ is a collection $\II$ of
$S_n$-invariant subspaces $\II(n) \subset \OO(n)$, for each $n \ge 1$,
such that whenever $i \in \II$, $\gamma(\dots, i, \dots) \in \II$.
\end{df}

The intersection of an arbitrary number of ideals in an operad is also
an ideal, and one can define the ideal generated by a subset in $\OO$
as the minimal ideal containing the subset.

\begin{df}
For an operad ideal $\II \subset \OO$, the \emph{quotient operad}
$\OO/\II$ is the collection $\OO(n)/\II(n)$, $n \ge 1$, with the
structure of operad induced by that on $\OO$.
\end{df}

The \emph{free operad $F(S)$ generated by a collection $S = \{ S(n) \;
  | \; n \ge 1 \}$ of sets}, is defined as follows.
\[
F(S)(n) = \bigoplus_{\text{$n$-trees $T$}} k \cdot S(T),
\]
where the summation runs over all planar rooted trees $T$ with $n$
labeled leaves and
\[
S(T) = \Map(v(T), S),
\]
the set of maps from the set $v(T)$ of vertices of the tree $T$ to the
collection $S$ assigning to a vertex $v$ with $\In(v)$ incoming edges
an element of $S(\In(v))$ (the edges are directed toward the root). In
other words, an element of $F(S)(n)$ is a linear combination of planar
$n$-trees whose vertices are decorated with elements of $S$. There is
a special tree with no vertices:

\centerline{\epsfysize=0.5in \epsfbox{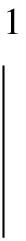}}
\smallskip

\noindent
The component $F(S)(1)$ contains, apart from $S(1)$, the
one-dimensional subspace spanned by this tree.

The following data defines an operad structure on $F(S)$.

\begin{enumerate}

\item The identity element is the special tree in $F(S)(1)$ with no 
vertices.

\item The symmetric group $S_n$ acts on $F(S)(n)$ by relabeling the inputs.

\item The operad composition is given by grafting the roots of trees
to the leaves of another tree. No new vertices are created.

\end{enumerate}

\begin{df}
  Now let $R$ be a subset of $F(S)$, \emph{i.e}., a collection of
  subsets $R(n) \subset F(S)(n)$. Let $(R)$ be the ideal in $F(S)$
  generated by $R$. The quotient operad $F(S)/(R)$ is called the
  \emph{operad with generators $S$ and defining relations $R$}.
\end{df}

\begin{ex}
  The associative operad $\ass$ is the operad generated by a point $S
  = S(2) = \{\bullet\}$ with a defining relation given by the
  associativity condition, see Section~\ref{ass:operad}, expressed in
  terms of trees.
\end{ex}

\begin{ex}
  The Lie operad $\lie$ is the operad also generated by a point $S =
  S(2) = \{\bullet\}$ with defining relations given by the skew symmetry
  and the Jacobi identity, see Section~\ref{lie:operad}.
\end{ex}

\begin{ex}
  The Poisson operad is the operad also generated by a two-point set
  $S = S(2) = \{\bullet, \circ\}$ with defining relations given by the
  commutativity and the associativity for simple trees decorated only
  with $\bullet$'s, the skew symmetry and the Jacobi identity for
  simple trees decorated with $\circ$'s, and the Leibniz identity for
  binary 3-trees with mixed decorations, see
  Section~\ref{poisson:operad}.
\end{ex}

\subsection{Homotopy algebra}

The idea of a \emph{homotopy ``something'' algebra} is to relax the
axioms of the \emph{``something'' algebra}, so that the usual
identities are satisfied up to homotopy. For example in a homotopy
associative algebra, the associativity identity looks like
\[
(ab)c - a(bc) \text{ is homotopic to zero}.
\]
Or in a homotopy Gerstenhaber (G-) algebra, the Leibniz rule is
\[
[a,bc] - [a,b]c \mp b[a,c] \text{ is homotopic to zero}.
\]
Usually, a homotopy something algebra arises when one wants to lift
the structure of a something algebra \emph{a priori} defined on \coh\
to the level of cochains.

This kind of relaxation seems to be too lax for many, practical and
categorical, purposes, and one usually requires that the
null-homotopies, regarded as new operations, satisfy their own
identities, up to their own homotopy. These homotopies should also
satisfy certain identities up to homotopy and so on. This resembles
Hilbert's chains of syzygies in early homological algebra. Algebras
with such chains of homotopies are called \emph{strongly homotopy
``something'' algebras} or \emph{``something''$_\infty$-algebras}.

Operads are especially helpful when one needs to work with
some\-thing$_\infty$-al\-ge\-bras. We already know that defining the
class of something algebras is equivalent to defining the something
operad. Thus, if we have an operad $\OO$, what is $\OO_\infty$,
\emph{the} corresponding strongly homotopy operad? M.~Markl's paper
\cite{markl} provides a satisfactory answer to this question: the
\emph{operad} $\OO_\infty$ is a minimal model of the operad $\OO$. A
minimal model is unique up to isomorphism. The idea is borrowed from
Sullivan's rational homotopy theory; a minimal model is, first of all,
a free resolution of $\OO$ in the category of operads of complexes,
\emph{i.e}., an operad of complexes free as an operad of graded vector
spaces, whose \coh\ is $\OO[0]$, the operad $\OO$ sitting in degree
zero, if it is an operad of vector spaces, or the operad $\OO$ sitting
in the original degrees, if it is already an operad of graded vector
spaces. Second of all, a minimal model must satisfy a minimality
condition: its differential must be decomposable.

For certain specific classes of operads, one manages to describe a
minimal model explicitly. For example, V.~Ginzburg and M.~Kapranov \cite{gk}
do it (even earlier than the notion of a minimal model for an operad
surfaced) for the so-called Koszul operads. Below we describe an
example of such kind, giving rise to the notion of an
$A_\infty$-algebra and the $A_\infty$ operad.

\subsubsection{$A_\infty$-algebras}

\begin{df}
  An {\it $A_\infty$-algebra}, or a {\it strongly homotopy associative
    algebra}, is a complex $V = \bigoplus_{i \in \nz} V^i$ with a
  differential $d$, $d^2 = 0$, of degree 1 and a collection of $n$-ary
  operations, called \emph{products}:
\[
M_n (v_1, \dots, v_n) \in V, \qquad v_1, \dots, v_n \in V,\; n \ge 2,
\]
which are homogeneous of degree $2-n$ and satisfy the relations
\begin{multline}
\label{Q}
d M_n(v_1, \dots, v_n) - (-1)^n \sum_{i=1}^n  \epsilon(i) M_n (v_1,
\dots, d v_i, \dots, v_n)
\\
= \sum_{\substack{k+l = n+1 \\ k, l \ge 2}}
\sum_{i=0}^{k-1} (-1)^{i+l(n-i-l)}
\sigma (i) M_k (v_{1}, \dots, v_{i},
M_l (v_{i+1}, \dots, v_{i+l}), \dots , v_n),
\end{multline}
where $\epsilon (i) = (-1)^{|v_1| + \dotsb + |v_{i-1}|}$ is the sign
picked up by taking $d$ through $v_1, \linebreak[0] \dots,
\linebreak[1] v_{i-1}$, $|v|$ denoting the degree of $v \in V$, and
$\sigma (i)$ is the sign picked up by taking $M_l$ through $v_1,
\linebreak[0] \dots, \linebreak[1] v_{i}$.
\end{df}

It is remarkable to look at these relations for $n =2$ and 3:
\[
d M_2(v_1,v_2) - M_2 (dv_1, v_2) - (-1)^{\abs{v_1}} M_2 (v_1, d v_2) = 
0,
\]
\begin{multline*}
d M_3(v_1,v_2,v_3) + M_3 (dv_1, v_2, v_3) +
(-1)^{\abs{v_1}} M_3 (v_1, d v_2, v_3)\\
\quad + (-1)^{\abs{v_1}+\abs{v_2}} M_3 (v_1, v_2, d v_3)\\
= M_2(M_2(v_1,v_2),v_3) - M_2(v_1, M_2(v_2,v_3)),
\end{multline*}
which mean that the differential $d$ is a derivation of the bilinear
product $M_2$ and the trilinear product $M_3$ is a homotopy for the
associativity of $M_2$, respectively.

$A_\infty$-algebras can be described as algebras over a certain tree
operad. This operad is the tree part of the graph complex, which will
be the topic of the following sections.

\subsubsection{The $A_\infty$ operad}

Let $A_\infty (n)$ be the linear span of the set of equivalence
classes of connected planar trees that have a root edge and $n$ leaves
labeled by integers 1 through $n$, with vertices of a valence at least
3, $n \ge 2$. For $n=1$ take one tree with a unique edge connecting a
leaf and a root.  Let us not include anything for $n=0$, although one
could do that similar to the associative operad case, so that the
corresponding notion of an $A_\infty$-algebra would have a unit.

We grade each vector space $A_\infty (n)$ by defining the degree $|T|$
of a tree $T \in A_\infty (n)$ via
\[
|T| := v(T) + 1 - n = e(T) + 1 - 2n,
\]
where $v(T)$ is the number of vertices and $e(T)$ the number of edges
of $T$. Notice that $2-n \le | T| \le 0$ for $n \ge 1$.

Let us define an operad structure on these spaces of trees. The
symmetric group acts by relabeling the leaves, and the operad
composition is obtained by \emph{grafting}, as in the examples above,
except one needs to take a sign into account. When we graft a tree
$T_2$ to the $i$th leaf of a tree $T_1$, the result must be the
grafted tree multiplied by a sign, which is $(-1)$ to the power
$(e(T_2)-1) ($the number of edges to the right of the $i$th leaf in
$T_1)$, where the edges to the right of a leaf are the edges which are
strictly on the right-hand side of a unique path from the leaf to the
root. The reason for the sign above is that grafting must respect the
differential, which is introduced below.

\begin{xca}
  Show that this operad is a free operad of vector spaces generated by
  the following trees for $n \ge 2$, which are sometimes called
  \emph{corollas}.

\begin{center}
\begin{picture}(180,65)(0,-12)
\put(-5,20){$\delta_n :=$}
\put(37,40){1}
\put(65,40){2}
\put(92,40){3}
\put(107,31){\dots}
\put(156,40){$n$}
\put(50,33){\vector(1,-1){27}}
\put(50,33){\line(1,-1){33}}
\put(72,33){\vector(1,-3){9}}
\put(72,33){\line(1,-3){11}}
\put(94,33){\vector(-1,-3){9}}
\put(94,33){\line(-1,-3){11}}
\put(147,32){\vector(-2,-1){53}}
\put(147,32){\line(-2,-1){64}}
\put(83,0){\vector(0,-1){12}}
\put(170,20){,}
\end{picture}
\end{center}

\end{xca}

\begin{rem}
  There is no need to mark directions on the edges of a tree: from now
  on we will assume the edges are directed from top to bottom.
\end{rem}

\subsubsection{The tree complex}

The above operad of trees is not yet the $A_\infty$-operad, but only
its underlying operad of graded vector spaces. The $A-\infty$-operad
is a DG operad, \emph{i.e}., an operad of complexes. The DG structure,
or a differential, is defined as follows.

Before defining it, we will define the operation of internal-edge
contraction on the set of trees.

\begin{df}
  We use the notation $T/e$ to denote the tree obtained from a tree
  $T$ by contracting an internal edge $e$:

\begin{center}
\begin{picture}(170,90)(0,-28)
\put(7,54){1}
\put(45,54){2}
\put(79,54){3}
\put(113,54){4}
\put(13,51){\line(1,-1){17}}
\put(30,34){\line(1,-1){17}}
\put(47,17){\line(1,-1){17}}
\put(47,51){\line(-1,-1){17}}
\put(81,51){\line(-1,-1){34}}
\put(115,51){\line(-1,-1){51}}
\put(64,0){\line(0,-1){12}}
\put(30,19){$e$}
\put(60,-28){$T$}
\end{picture}
\hspace{\bigskipamount}
\begin{picture}(170,90)(0,-28)
\put(7,54){1}
\put(45,54){2}
\put(79,54){3}
\put(113,54){4}
\put(13,51){\line(1,-1){34}}
\put(47,17){\line(1,-1){17}}
\put(47,51){\line(0,-1){34}}
\put(81,51){\line(-1,-1){34}}
\put(115,51){\line(-1,-1){51}}
\put(64,0){\line(0,-1){12}}
\put(60,-28){$T/e$}
\end{picture}
\end{center}
\end{df}

\noindent

We can now define a {\it differential} $d: A_\infty (n) \to A_\infty
(n)$ by the formula
\[
dT := \sum_{T': T = T'/e} \epsilon T' ,
\]
where $\epsilon$ is the sign given by counting the number of edges
below and to the left of the edge $e$ in the tree $T'$, not counting
the root.

In particular,
\begin{multline}
\label{d}
\quad
d \quad
\begin{picture}(119,60)(37,8)
\put(37,40){1}
\put(65,40){2}
\put(92,40){3}
\put(107,31){\dots}
\put(156,40){$n$}
\put(83,0){\line(-1,1){33}}
\put(83,0){\line(-1,3){11}}
\put(83,0){\line(1,3){11}}
\put(83,0){\line(2,1){63}}
\put(83,0){\line(0,-1){12}}
\end{picture}
\\
=
\sum_{  \substack{
   k+l = n+1\\
   k, l \ge 2}}
\sum_{i=0}^{k-1} \; (-1)^i 
\begin{picture}(119,60)(37,8)
\put(43,40){1}
\put(60,40){$i+1$}
\put(98,40){$i+l$}
\put(151,40){$n$}
\put(83,0){\line(-1,1){33}}
\put(55,31){\dots}
\put(83,0){\line(1,3){11}}
\put(117,31){\dots}
\put(90,21){\line(-1,1){12}}
\put(90,21){\line(-1,3){4}}
\put(90,21){\line(3,2){18}}
\put(83,0){\line(1,1){33}}
\put(83,0){\line(2,1){66}}
\put(83,0){\line(0,-1){12}}
\end{picture} .
\end{multline}

\begin{prop}
\begin{enumerate}
\item The operator $d$ satisfies $d^2 = 0$ and $\deg d = 1$.
\item The operad structure on $A_\infty = \{ A_\infty (n) \; | \; n
  \ge 1\}$ is compatible with the differential $d$:
\[
d(T_1 \circ_i T_2) = d T_1 \circ_i T_2 + (-1)^{\abs{T_1}} T_1 \circ_i
dT_2,
\]
\emph{i.e.,} $A_\infty$ is a DG operad.
\end{enumerate}
\end{prop}

\begin{df}
  We will call the DG operad $A_\infty$ the {\it $A_\infty$ operad}.
\end{df}

\begin{rem}
\label{Lie}
The complex $A_\infty (n)$ is part of the (cochain) graph complex, see
Section~\ref{graph-homology}. A similar operad $L_\infty$, based on
abstract, \emph{i.e}., nonplanar trees, was introduced by V.~Hinich
and V.~Schechtman \cite{hs}.  The operad $A_\infty$ is the dual cobar
operad in the sense of Ginzburg and Kapranov \cite{gk} of the
associative operad $\ass$. They also show that the \coh\ of the operad
$A_\infty$ is the associative operad $\ass$ of
Section~\ref{ass:operad}, implying that $A_\infty$ is a free, and in
fact, minimal, resolution of $\ass$.
\end{rem}

The following theorem shows that the $A_\infty$ operad describes the
class of $A_\infty$-algebras.
\begin{thm}[\cite{gk}]
%\label{hass}
  An algebra over the $A_\infty$ operad is an $A_\infty$-algebra. Each
  $A_\infty$-algebra admits a natural structure of an algebra over the
  $A_\infty$ operad.
\end{thm}

\begin{proof}
  For a complex $V$ of vector spaces with a differential $d$ of degree
  1, $d^2 =0$, the structure of an algebra over the operad $A_\infty$
  on $V$ is a morphism of DG operads:
\[
\phi: A_\infty (n) \to \End{V}(n), \qquad n \ge 1,
\]
where $\End{V}(n) := \Hom (V^{\otimes n}, V)$ is the {\it endomorphism
  operad}, which is also a DG operad (with the usual internal
differential determined by $d$). Given such a morphism $\phi$, we
define the $n$-ary product on $V$:
\[
M_n(v_1, \dots , v_n) := \phi(\delta_n)(v_1 \otimes \dotsb \otimes
v_n).
\]
Note that the degree of the product is equal to that of the corolla
$\delta_n$, which is $2 - n$.  Since $\phi$ is a morphism of DG
operads, $d\phi = \phi d$, and in view of \eqref{d}, this is
equivalent to the identity \eqref{Q}.

Conversely, given a collection of $n$-ary brackets on $V$, $n \ge 2$,
we define a morphism $\phi$ on the generators $\delta_n$ by the above
formula. The $A_\infty$ operad is freely generated by the corollas
$\delta_n$, with a differential defined by \eqref{d}, so the mappings
$\phi$ define a morphism of DG operads, if the relations \eqref{d} are
satisfied by the $\phi(\delta_n)$'s. Equations \eqref{Q} show that
this is the case.
\end{proof}

\subsection{Cacti}
\label{cacti}

One of the very recent applications of operads moves somewhat
southwest and goes through cacti rather than trees. However, the word
``homology'' will mean rational homology throughout this section.

The cacti operad may be used to explain the structure of a BV-algebra
on the homology of a free loop space on a compact oriented manifold
discovered by M.~Chas and D.~Sullivan \cite{chas-sullivan}. We will be
brief and not describe what exactly a BV-algebra structure is.  Even
if you are new to the subject, you may bear with me by taking a
tautological approach: think of BV as the algebraic structure induced
on homology from the structure of an algebra over the cacti operad on
a topological space. Some of the missing details are also captured by
Ralph Cohen's lecture and his paper with John Jones
\cite{cohen-jones}, where they develop a more civilized version of the
cacti operad action in the category of spectra and show that
Chas-Sullivan's BV-structure is the same as the BV-structure coming
naturally after identifying the homology of the free loop space with
the Hochschild homology of the cochain algebra of the target space.

The \emph{cacti operad} is an operad $c = \{ C(n) \; | \; n \ge 1\}$
of topological spaces. Its $n$th component $C(n)$ for $n \ge 1$ may be
described as follows.

$C(n)$ is the set of ordered tree-like configurations of parameterized
circles (lobes) of varying (positive) radii, along with a cyclic order
of components at the intersection points and the choice of a point
``0'' on the whole configuration along with the choice of one of the
circles on which this point 0 lies. The latter is, of course essential
when 0 happens to be an intersection point. The topology on the set of
configurations before choosing cyclic orders and marking a point 0 is
induced from a natural embedding into $(S^1)^{{n \choose 2}}$. The
choice of cyclic orders defines a finite covering. After these choices
are made, we can define a continuous flow on the cactus which goes
along the parameters on the circles and jumps from one component to
the next in the cyclic order at the intersection points. Two choices
of a zero point are considered close to each other, if they are close
in the sense of this flow. Thus the space of cacti with a marked
point 0 is an $S^1$-bundle over the configuration space. In short,
before marking 0 the topology we are talking about is Gromov's
topology: ``two pictures are close to each other, if they appear
so''. This principle is quite subtle and does not always work the same
way: a classical example of its violation is perhaps the right side
mirror on your car, cf.\ \cite{jurassic}.
\smallskip

\centerline{\epsfxsize=1.5in \epsfbox{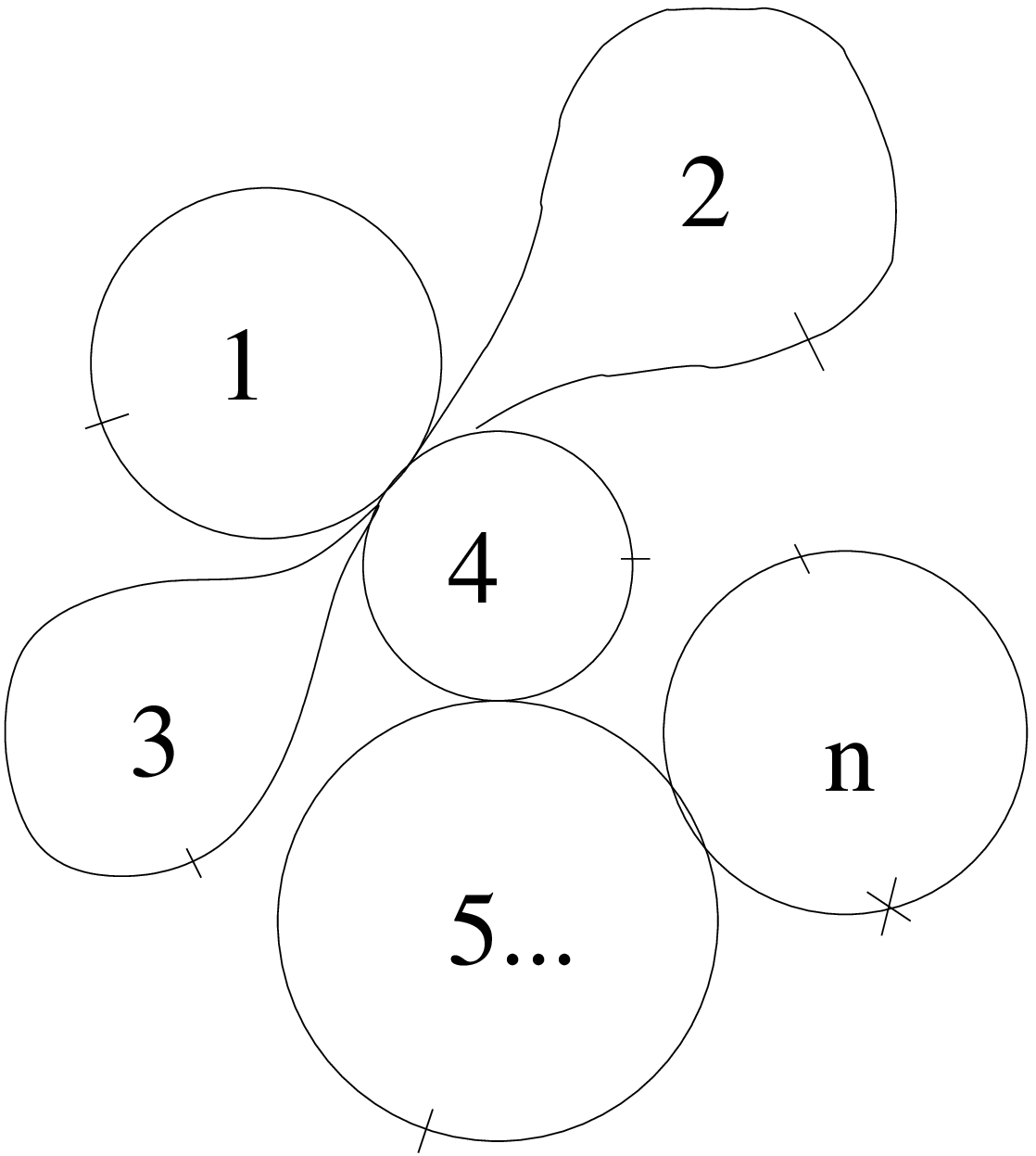}}
\smallskip

The operad structure on the cacti comes from the following
observation. The choice of a point 0 and a component on a cactus gives
a natural map from $S^1$ to the cactus. First rescale the radius of
$S^1$ to match the sum of the radii of the lobes forming the
cactus. Then wind this $S^1$ around the cactus and follow the flow
along the lobes, starting with the chosen lobe at 0. Topologically,
the constructed map will identify a few groups of points on $S^1$ and
therefore will have a degree one on each lobe. Given two cacti and the
$i$th lobe in the first one, the operad composition $\circ_i$ will be
given by further collapsing the $i$th circle according to the map
given by the second cactus.

The following theorem describes both the cacti operad $C$ and how it
produces the BV-structure on the homology of a free loop space $LM$ in
a compact oriented manifold $M$ of dimension $d$.

\begin{thm}
 \begin{enumerate}

 \item
The cacti operad is homotopy equivalent to the framed little disks
operad.

 \item
The cacti operad $C$ ``acts'' on $LM$ in the following sense. The
diagram
\[
C(n) \times (LM)^n \xleftarrow{i} L^{C(n)} M \xrightarrow{e} LM,
\]
where $L^{C(n)} M$ is the space of continuous maps of $n$-component
cacti to $M$, induces a composite map
\[
H_\bullet (C(n)) \otimes H_\bullet(LM)^{\otimes n} \xrightarrow{i^!}
(H_\bullet (L^{C(n)} M))[(1-n)d] \xrightarrow{e_*} H_\bullet
(LM)[(1-n)d],
\]
where $i^!$ denotes the pullback in homology. The collection of such
maps for $n \ge 1$ is compatible with the operad structure on $C$.

\item
The composite map $e_* i^!$ produces the structure of an algebra over
the homology cacti operad $H_\bullet (C)$ on the space
$H_\bullet(LM)[d]$.

\end{enumerate}
\end{thm}

Combining Statements 1 and 3 of the theorem with a theorem of
E.~Getzler \cite{g} which says that the homology framed little disks
operad is the operad describing BV-algebras and checking what the
basic operations (the dot product and the BV operator) really are, we
obtain the following result.

\begin{crl}
The space $H_\bullet(LM)$ (after an appropriate degree shift) has the
natural structure of a BV-algebra, coinciding with the one constructed
by Chas and Sullivan.
\end{crl}

Cacti of higher dimensions are much subtler and are the topic of an
upcoming paper with Sullivan. Here we will just mention one of the
applications, which uses a generalization of Getzler's theorem to
higher dimensions by P.~Salvatore and N.~Wahl \cite{salvatore-wahl}.

\begin{thm}[Sullivan-A.V. \cite{sullivan-me}]
Let $M$ be a compact oriented manifold and $S^n M = \Map(S^n, M)$ the
sphere space. Then the space $H_\bullet (S^n M)$ (after an appropriate
degree shift) is naturally an algebra over the homology framed little
$(n+1)$-disks operad. In particular, the homology of the sphere space
has the following algebraic structure. For $n$ odd, it is a
BV$_{n+1}$-algebra, which is the same as a usual BV-algebra, except
that the BV-operator has degree $n$. For $n$ even, this structure is
the same as that of an $(n+1)$-algebra with a differential of degree
$2n-1$.
\end{thm}

\section{Graph homology}
\label{graph-homology}

Even more interesting things start happening when you pass from trees
to graphs. On of these things is that such a classical, analytic and
algebraic geometric object as the moduli space of Riemann surfaces
miraculously emerges in the horizon the very moment you say the word
``graphs''. Out of the numerous versions of the graph complex, we have
chosen the ribbon one, which is most closely related to the moduli
spaces of Riemann surfaces. The notion of a ribbon graph and the graph
complex are due to R.~Penner \cite{penner:bull,penner:fatgraph}, who
used the term ``fatgraph''. The ribbon graph complex is a
generalization of the planar tree complex we considered in the
previous section, except that we no longer allow any free legs. Other
versions of the graph complex include those which do not require
cyclic orders at vertices, see M.~Culler and K.~Vogtmann
\cite{culler-vogtmann}, or, on the contrary, have more complicated
decorations at vertices, see Getzler-Kapranov's Feynman transform
\cite{gek}. There is also a dual version of the graph complex,
producing graph \coh. Kontsevich \cite{kon:sympl}, who noticed a
common pattern in Penner's and Culler-Vogtmann's work and related
different versions of the graph complex to different fundamental types
of algebras (or operads that describe them), had an enormous influence
on the subject.

\subsection{The graph complex}

The right analogue of a planar tree is a \emph{ribbon graph}, which is
a nonempty connected finite cell complex $\Gamma$ of dimension one
with the choice of a cyclic order on the set of half-edges around each
vertex. We also require that the valences of vertices must be greater
than or equal to 3. A \emph{boundary component} of a ribbon graph
$\Gamma$ is a cyclic sequence $\vec e_0, \vec e_1, \dots, \vec e_q =
\vec e_0$ of directed edges (\emph{i.e}., edges with directions chosen
on each of them) of $\Gamma$, so that for each pair $(\vec e_i, \vec
e_{i+1})$ of two subsequent edges, the tail of $\vec e_{i+1}$ is the
half-edge that follows the head of $\vec e_i$ in the cyclic order at a
common vertex of $e_i$ and $e_{i+1}$. To have a more direct connection
with moduli spaces of Riemann surfaces with labeled boundary
components, we will assume that graphs under consideration will have
their boundary components labeled. We identify ribbon graphs which
differ by an \emph{isomorphism}, which is just a cellular
homeomorphism preserving the cyclic orders at vertices and fixing each
of the boundary components (but not necessarily each edge forming a
boundary component). The \emph{group of automorphisms} of a ribbon
graph is denoted by $\Aut_\del (\Gamma)$. This group is infinite, but
has a finite number of connected components, $\pi_0 (\Aut_\del
(\Gamma))$, which is what would usually be called the group of
automorphisms of a graph. An \emph{orientation} on a graph $\Gamma$ is
an orientation on the vector space $\nr^{e(\Gamma)}$, where
$e(\Gamma)$ is the set of edges of $\Gamma$.

For each $m \ge 2$, let $G_m$ be the free abelian group generated by
the set of isomorphism classes of oriented ribbon graphs $\Gamma$ as
above with $m$ edges, modded out by the defining relations $\Gamma
+ (-\Gamma) = 0$, where $-\Gamma$ is the same graph as $\Gamma$, but
with the opposite orientation.

Define a differential $d: G_m \to G_{m-1}$, so that $d^2 =0$, as
follows:
\[
d \, \Gamma := \sum_{\substack{\text{edges $e \in e(\Gamma)$}\\
\text{which are not loops}}} \Gamma/e,
\]
where $\Gamma/e$ is the ribbon graph obtained from $\Gamma$ by
contracting edge $e$ to a point. The cyclic order at the new vertex
created by merging the two ends of $e$ is obtained by a natural
insertion. The orientation is defined by the natural isomorphism
\begin{eqnarray*}
\Lambda^{\max} \nr^{e(\Gamma/e)} & \to & \Lambda^{\max} \nr^{e(\Gamma)},\\
\omega & \mapsto &  \omega \wedge e .
\end{eqnarray*}
The complex $G_\bullet$ with the differential $d$ is called the
\emph{graph $($chain$)$ complex} and its homology is called
\emph{graph homology}.

Obviously, the differential preserves the number $n$ of boundary
components, as well as the \emph{Euler characteristic} $\chi (\Gamma)
:= v(\Gamma) - e (\Gamma)$. It also preserves the \emph{genus} $g$
defined by the equation $\chi = 2 - 2g -n$. The solution $g= 1 - (\chi
+ n)/2$ is a nonnegative integer, because if we glue in $n$ disks into
the boundary components, we will get a compact orientable topological
surface whose Euler characteristic equals $\chi + n = 2-2g$ for some
nonnegative integer $g$. Thus, the graph complex splits into the direct
sum $G_\bullet = G_\bullet^{g,n}$ of subcomplexes $G_\bullet ^{g,n}$
with a fixed genus $g$ and a number $n$ of boundary components, $g \ge
0$, $n \ge 1$.

\subsection{Metric ribbon graphs and moduli spaces}
  
The graph complex above is in fact the chain complex of a certain cell
complex, that of metric ribbon graphs. A \emph{metric} on a ribbon
graph is an assignment of a positive real number, a \emph{length}, to
each edge of the graph. Obviously, the space of isomorphism classes of
metric ribbon graphs with an underlying graph $\Gamma$ is the space
$\nr_+^{e(\Gamma)}/ \Aut_\del (\Gamma)= \nr_+^{e(\Gamma)}/
\pi_0(\Aut_\del (\Gamma))$, which is, in fact, an orbifold. For
different ribbon graphs $\Gamma$ of a fixed genus $g$ and a number $n$
of boundary components, these orbifolds glue together by identifying
metric graphs some of whose lengths degenerate to zero with the metric
graphs obtained by contracting the zero length edges.  Thus, these
orbifolds become what is called the \emph{rational cells} of an
orbifold, the \emph{space} $G^{\met}_{g,n}$ \emph{of metric ribbon
  graphs}.  If we forget the orbifold structure, these rational cells
are just cells of a topological space, but not a cell complex, because
only part of the boundary of each cell, the part with lengths
approaching zero, is glued up to other cells. To get an honest cell
complex, we may take the one-point compactification of the space of
metric ribbon graphs, thus, gluing in a single point to all cells as
lengths tend to positive infinity. This gives a nonorbifold-type
singularity, and we get a compact smooth orbifold $\overline
G^{\met}_{g,n}$ with one singular point.  It may also be viewed as an
ordinary cell complex with a base point.  This space was designed with
the following obvious proposition in mind.

\begin{prop}
\label{GM}
The chain complex computing the rational reduced homology of the
orbifold $\overline G^{\met}_{g,n}$ is isomorphic to the ribbon graph
complex $G_\bullet^{g,n} \otimes \nq$.
\end{prop}

\begin{rem}
  We had to take the rational coefficients to make sure that we were
  computing the orbifold homology. Over the integers, it would be the
  computation of the ordinary homology of the underlying space, which
  would be somewhat misleading.
\end{rem}

The relevance of the graph complex and the metric ribbon graph space
is explained by the following theorem, which is one of the deepest
results in mathematics in the past twenty years.

\begin{thm}[Harer-Mumford-Thurston, Penner, see a remark below]
  For $g, n > 0$ (or $g = 0$ and $n > 2$), there is an orbifold
  isomorphism
\[
\MM_{g,n} \times \nr^n_+ \cong G^{\met}_{g,n},
\]
where $\MM_{g,n}$ is the moduli space of compact smooth Riemann
surfaces of genus $g$ with $n$ labeled punctures.
\end{thm}

\begin{rem}
  The idea of such combinatorial description of the moduli space
  belongs to D.~Mumford and W.~Thurston. Mumford's approach, which was
  realized by J.~Harer \cite{harer:vir-dim}, uses the theory of
  Strebel differentials on a Riemann surface, an analytic result
  producing a unique meromorphic quadratic differential $q$ in such a
  way that a metrized ribbon emerges as the set of critical horizontal
  trajectories of $q$ and the residues of $\sqrt{q}$ at punctures give
  elements in $\nr^n_+$, see more detail in E.~Looijenga
  \cite{looijenga:texel}, R.~Hain-Looijenga \cite{hain-looijenga}, and
  M.~Mulase-M.~Penkava \cite{mulase-penkava}. Penner \cite{penner:CMP}
  came up with a quite different way to obtain the above result, in
  which the points of $\nr^n_+$ on the left-hand side were interpreted
  as the hyperbolic lengths of horocycles attached to the punctures in
  a hyperbolic model of the Riemann surface with the punctures
  removed. Kontsevich \cite{kon:witten} considered a natural
  compactification of the space of metric ribbon graphs and compared
  it with Deligne-Mumford's one for $\MM_{g,n} \times \nr^n_+$.
\end{rem}

Combining this theorem with Proposition~\ref{GM}, one gets the
following purely combinatorial way of computing the rational homology
of the moduli spaces.

\begin{crl}
  The homology of the ribbon graph complex $G_\bullet^{g,n} \otimes
  \nq$ is isomorphic to the rational reduced homology of the one-point
  compactification of $\MM_{g,n} \times \nr_+^n$.
\end{crl}

%\begin{sloppypar}
Applying Poincar\'{e}-Lefschetz duality, we see that
\[
H_k(G_\bullet^{g,n} \otimes \nq) \linebreak[1] = \linebreak[0]
H^{6g-6+3n-k} (\MM_{g,n}; \nq).
\]
%\end{sloppypar}

The above corollary has been used to prove Harer's stability theorem,
estimate the homological dimension of $\MM_{g,n}$, and compute its
virtual (\emph{i.e}., orbifold) Euler characteristic, see
Hain-Looijenga \cite{hain-looijenga}.  Unfortunately, the
combinatorics of the graph complex is complicated enough not to
produce more general, explicit computations of its homology as of yet.
However, this combinatorial model was one of the key elements in
Kontsevich's proof \cite{kon:witten} of the Witten conjecture
\cite{witten} on the intersection theory on moduli spaces and the KP
hierarchy.

\bibliographystyle{amsalpha}
\bibliography{op}
\end{document}